\documentclass{article}
\usepackage{array,amsmath}
\numberwithin{equation}{section}
\usepackage{indentfirst}

\begin{document}

\newtheorem{thm}{Theorem}[section]
\newtheorem{cor}[thm]{Corollary}
\newtheorem{prop}[thm]{Proposition}
\newtheorem{conj}[thm]{Conjecture}
\newtheorem{lem}[thm]{Lemma}
\newtheorem{Def}[thm]{Definition}
\newtheorem{rem}[thm]{Remark}
\newtheorem{prob}[thm]{Problem}
\newtheorem{ex}{Example}[section]

\newcommand{\be}{\begin{equation}}
\newcommand{\ee}{\end{equation}}
\newcommand{\ben}{\begin{enumerate}}
\newcommand{\een}{\end{enumerate}}
\newcommand{\beq}{\begin{eqnarray}}
\newcommand{\eeq}{\end{eqnarray}}
\newcommand{\beqn}{\begin{eqnarray*}}
\newcommand{\eeqn}{\end{eqnarray*}}
\newcommand{\bei}{\begin{itemize}}
\newcommand{\eei}{\end{itemize}}

\newcommand{\pa}{{\partial}}
\newcommand{\V}{{\rm V}}
\newcommand{\R}{{\bf R}}
\newcommand{\K}{{\rm K}}
\newcommand{\e}{{\epsilon}}
\newcommand{\tomega}{\tilde{\omega}}
\newcommand{\tOmega}{\tilde{Omega}}
\newcommand{\tR}{\tilde{R}}
\newcommand{\tB}{\tilde{B}}
\newcommand{\tGamma}{\tilde{\Gamma}}
\newcommand{\fa}{f_{\alpha}}
\newcommand{\fb}{f_{\beta}}
\newcommand{\faa}{f_{\alpha\alpha}}
\newcommand{\faaa}{f_{\alpha\alpha\alpha}}
\newcommand{\fab}{f_{\alpha\beta}}
\newcommand{\fabb}{f_{\alpha\beta\beta}}
\newcommand{\fbb}{f_{\beta\beta}}
\newcommand{\fbbb}{f_{\beta\beta\beta}}
\newcommand{\faab}{f_{\alpha\alpha\beta}}

\newcommand{\pxi}{ {\pa \over \pa x^i}}
\newcommand{\pxj}{ {\pa \over \pa x^j}}
\newcommand{\pxk}{ {\pa \over \pa x^k}}
\newcommand{\pyi}{ {\pa \over \pa y^i}}
\newcommand{\pyj}{ {\pa \over \pa y^j}}
\newcommand{\pyk}{ {\pa \over \pa y^k}}
\newcommand{\dxi}{{\delta \over \delta x^i}}
\newcommand{\dxj}{{\delta \over \delta x^j}}
\newcommand{\dxk}{{\delta \over \delta x^k}}

\newcommand{\px}{{\pa \over \pa x}}
\newcommand{\py}{{\pa \over \pa y}}
\newcommand{\pt}{{\pa \over \pa t}}
\newcommand{\ps}{{\pa \over \pa s}}
\newcommand{\pvi}{{\pa \over \pa v^i}}
\newcommand{\ty}{\tilde{y}}
\newcommand{\bGamma}{\bar{\Gamma}}

\font\BBb=msbm10 at 12pt
\newcommand{\Bbb}[1]{\mbox{\BBb #1}}

\newcommand{\qed}{\hspace*{\fill}Q.E.D.}  

\title{ The navigation problems and the curvature properties on conic Kropina manifolds}
\author{ Xinyue Cheng \footnote{supported by the National Natural Science Foundation of China (11871126) and the Science Foundation of Chongqing Normal University (17XLB022)}, Qiuhong Qu, Suiyun Xu}

\maketitle

\begin{abstract}
In this paper, we study navigation problems on conic Kropina manifolds. Let $F(x, y)$ be a conic Kropina metric on an $n$-dimensional manifold $M$ and  $V$ be a conformal vector field on  $(M, F)$ with $F(x, - V_{x})\leq 1$. Let $\widetilde{F}= \widetilde{F} (x,y)$ be the solution of the navigation problem with navigation data $(F, V)$.  We prove that $\widetilde{F}$ must be either a Randers metric or a Kropina metric. Then we establish the relationships between some curvature properties of $F$ and the corresponding properties of the new  metric $\widetilde{F}$, which involve S-curvature, flag curvature and Ricci curvature.
\end{abstract}

\section{Introduction}

In 1931, E. Zermelo studied the following problem (\cite{EZ}):

{\it Suppose that a ship sails on the open sea in calm waters and a mild breeze comes up. How must the ship be steered in order to reach a given destination in the shortest time?}

The problem was solved by Zermelo himself for the Euclidean flat plane and by Z. Shen (\cite{ZS2}\cite{ZS3}) in the case when the sea is a Riemannian manifold $(M, h)$ under the assumption that
the wind $W$ is a time-independent mild breeze, i.e. $h(x, W) < 1$.

Essentially, Zermelo navigation problem is tightly related to the geometry of indicatrix. Let $(M, \Phi)$ be a Finsler manifold. For each $x \in M$, the indicatrix  $S_{\Phi}$ of $\Phi$ at $x$ is a closed hypersurface of $T_{x}M$ around the origin $x$ defined by
\[
S_{\Phi}=\{y\in T_{x}M \ | \Phi (x,y)= 1 \}.
\]
Let $ W =W(x)$ be a  vector field on $M$. Consider the parallel shift $S_{\Phi}+\{ W \}$ of $S_{\Phi}$ along $W$.  It is easy to see that $y \in S_{\Phi}+\{ W \}$ if and only if $\Phi (x, y-W_{x})=1$. Further,  when $\Phi (x, -W_{x})<1$, $S_{\Phi}+\{ W \}$ contains the origin $x$ of $T_{x}M$. In this case, for any $y \in T_{x}M \setminus \{0\}$, there is a unique positive number $F=F(x,y) >0$ such that
\be
\frac{y}{F(x,y)}\in S_{\Phi}+\{ W \}, \label{naviga1}
\ee
that is, $F=F(x,y)$ satisfies the following
\be
\Phi \left( x, \frac{y}{F(x,y)}-W_{x}\right) = 1. \label{navigation}
\ee
It is easy to see that $F=F(x,y)$ is a regular Finsler metric (\cite{ChernShen}). On the other hand, when $\Phi (x, -W_{x})= 1$, $S_{\Phi}+\{ W \}$ passes thorough the origin $x$ of $T_{x}M$. In this case, for any $y\in A_{x}:=\left\{ y\in T_{x}M \ | \ g_{W}(y, W)>0 \right\}$, there is a unique positive number $F=F(x,y)$ such that $F=F(x,y)$ satisfies (\ref{naviga1}),  equivalently, $F= F(x,y)$ satisfies (\ref{navigation}). Actually, such $F=F(x,y)$ is a conic Finsler metric. The conic Finsler metric $F=F(x,y)$ is regular on $A:=\bigcup_{x\in M}A_{x} \subset TM$. In general, Finsler metric $F=F(x,y)$ obtained from (\ref{navigation}) is called a solution of the Zermelo navigation problem with navigation data $(\Phi , W)$.

A Finsler manifold $(M, F)$ of dimension $n$ is called  $C$-reducible if the following three conditions are satisfied: (1) \ $F$ is not Riemannian; (2) \ The dimension $n \geq 3$; (3) \ The Cartan torsion ${\bf C}$ of $F$ is written in the form
\[
C_{ijk}==\frac{1}{n+1}\left\{I_{i} h_{j k}+I_{j} h_{i k}+I_{k} h_{i j}\right\},
\]
where $h_{i j} :=F F_{y^{i} y^{j}}= g_{ij}-F_{y^{i}}F_{y^{j}}$ denote the angular metric tensor of $(M, F)$ and $I_{i}=g^{j k} C_{i j k}$ denote the mean Cartan torsion.
In this case, $F$ is called the  $C$-reducible Finsler metric . A fundamental fact is that any $C$-reducible Finsler metric is of the Randers type or the Kropina type only (see \cite{MM}). Randers metrics are one of the simplest non-Riemannian Finsler metrics with the form $F=\alpha+\beta$, where $\alpha=\sqrt{a_{ij}(x)y^{i}y^{j}}$ is a Riemannian metric and $\beta = b_{i}(x)y^{i}$ is a 1-form with $\|\beta _{x}\|_{\alpha}<1$ on the manifold. Randers metrics arise from the general relativity. Kropina metrics form a special kind of Finsler metrics in the form $F=\frac{\alpha ^2}{\beta}$. Kropina metrics were first introduced by L. Berwald when he studied the two-dimensional Finsler spaces with rectilinear extremal and investigated by Kropina (see \cite{VK1}\cite{VK2}). Kropina metrics have important and interesting applications in the theory of thermodynamics. Besides, both of Randers metrics and Kropina metrics play an interesting role in the Krivan problem in ecology (\cite{AIM}).  However,  Randers metrics are regular Finsler metrics, but Kropina metrics are the Finsler metrics with singularity.  In fact, Kropina metrics are not classical Finsler metrics, but conic Finsler metrics (\cite{XQ2}\cite{RS}).

Randers metrics and Kropina metrics can be both expressed as the solution of the Zermelo navigation problem on some Riemannian manifold $(M, h)$ with a vector field $W$. Concretely, assume that $h=\sqrt{h_{ij}(x)y^iy^j}$ and  $W=W^{i}\frac{\partial}{\partial x^{i}}$ with $\|W\|_h < 1$, then the metric $F$ obtained by solving the following problem
\be
h\left(x,\frac{y}{F(x,y)} -W_{x}\right) =1     \label{eqa2}
\ee
is a Randers metric given  by
\be
F=\frac{\sqrt{\lambda h^{2}+ W_{o}^{2}}}{\lambda}-\frac{W_{0}}{\lambda}. \label{RhW}
\ee
Here, $W_{0}:=W_{i}y^{i}=  h (y,W_{x})$, $W_{i}:=h_{ij}W^{j}$, $\lambda:=1-\|W\|_{h}^{2}>0$ and $\|W_{x}\|_{h}=\|\beta_{x}\|_{\alpha}$. The condition $\|W\|_h < 1$ is essential for obtaining a positive definite Randers metric by  Zermelo navigation problem. In this case, we call $(h, W)$ the navigation data of Randers metric $F=\alpha +\beta$ (see \cite{BRS}\cite{ChernShen}). On the other hand, assume that  $W=W^i \frac{\pa }{\pa x^i}$ is a vector field with $\|W\|_h =1$. Then the solution of the Zermelo navigation problem (\ref{eqa2}) is a Kropina metric given by
\be
F=\frac{h^2}{2W_0}. \label{eqa3}
\ee
In fact, $F$ given by (\ref{eqa3}) is a conic Kropina metric defined on the conic domain
\[
A =\big\{(x,y) \in TM \ | \ h (y,W_{x}) >0\big\}\subset TM.
\]
In this case, $(h, W)$ is called the navigation data of conic Kropina metric $F=\frac{\alpha ^2}{\beta}$ and
\be
\alpha^2=\frac{b^2}{4}h^2, \ \ \ \beta=\frac{b^2}{2}W_0.\label{eqa4}
\ee
Here, $b :=\| \beta \|_\alpha$ denotes the norm of $\beta$ with respect to $\alpha$.  Conversely, given a conic Kropina metric $F=\frac{\alpha^2}{\beta}$, put
\be
h_{ij}=\frac{4}{b^2}a_{ij},  \ \ \ W^{i}=\frac{1}{2}b^{i}. \label{NaRe}
\ee
Then we can get a Riemannian metric $h$ and a vector field $W$ with $\|W\|_{h}=1$ from (\ref{NaRe}) and $F$ is just given by (\ref{eqa2}) for $h$ and $W$.
Thus there is an one-to-one correspondence between a conic Kropina metric $F$ and a pair $(h,W)$ with $\|W\|_{h}=1$.  It is easy to see that a conic Kropina metric can be regarded as the limit of the navigation problem for Randers metrics when $\|W\|_h \to 1$. In the following, we just study conic Kropina metrics and we always use Kropina metric to take place of conic Kropina metric.

The theory on conformal vector fields is one of the core contents of conformal geometry, which plays an important role in differential geometry and physics.
In \cite{HMo}, Huang-Mo determined the flag curvature of a Finsler metric $\widetilde{F}$ produced from a Finsler metric $F$ and its  homothetic field $V$ in terms of the navigation problem. In \cite{SX}, Shen-Xia studied Zermelo navigation problem on a Randers manifold $(M, F=\alpha +\beta)$ with a conformal vector field $V$. The solution of the Zermelo navigation problem with navigation data $(F, V)$ is a new Randers metric $\widetilde{F}$. Then they established the relationships between some curvature properties of $F$ and the corresponding properties of $\widetilde{F}$, in which $S$-curvature, flag curvature and Ricci curvature were involoved. Further, Q. Xia considered a Randers metric $F$ of Douglas type on a manifold $M$  and a conformal vector field $V$ on $(M, F)$. Let $\widetilde{F}$ be the Randers metric generated from the navigation data $(F, V)$ by solving navigation problem. Q. Xia established the relationship among the flag curvature of $F$, the flag curvature of $\widetilde{F}$ and the conformal vector field $V$ (see \cite{XQ1}).

Motivated by above researches and based on the characterizations of conformal vector fields on Korpina manifolds given in \cite{CLY}, we study Zermelo navigation problem on a Kropina manifold $(M, F=\frac{\alpha ^2}{\beta})$. For a vector field $V$ on Kropina manifold $(M, F)$ with $-V_{x}\in A_{x}:=\{y\in T_{x}M \ | \ \beta = b_{i}(x) y^{i}>0\}$, let $\widetilde{F}$ be the solution of the Zermelo navigation problem with navigation data $(F, V)$. We firstly prove the following results: when $F(x, -V_{x})<1$,  $\widetilde{F}$ is a Randers metric; when $F(x, -V_{x})=1$,  $\widetilde{F}$ is a Kropina metric (for the details, see Theorem \ref{lem3}). Then, we can prove the following theorems.

\begin{thm}\label{Kropina navigation} \ Let $F=\frac{\alpha^2}{\beta}$ be a Kropina metric on an $n (\geq 2)$-dimensional manifold $M$ expressed by (\ref{eqa3}) with the navigation data $(h, W)$. Assume that $V$ is a conformal vector field on $(M, F)$ with conformal factor $\rho(x)$, $- V_{x}\in A_{x}$ and $F(x, -V_{x})<1$. Let $\widetilde{F}=\widetilde{F}(x, y)$ be the Randers metric defined from $(F, V)$ by
\be
F\left(x, \frac{y}{\widetilde{F}(x,y)}- V_{x}\right)=1, \ \ \ y \in T_{x}M. \label{KropinaNaviga}
\ee
Then we have the following.
\ben
\item[{\rm (1)}] If $F$ is of isotropic $S$-curvature with ${\bf S}(x, y)=(n+1)cF(x, y)$, then $\widetilde{F}$ is of isotropic $S$-curvature with $\widetilde{\bf S}(x, u)=(n+1)(- \rho)\widetilde{F}(x, u)$. In this case, $c$ is actually zero. Here and below, $u:=y+F(x, y)V=y+\widetilde{F}(x, u)V$.
\item[{\rm (2)}] If $F$ is of weakly isotropic flag curvature with ${\bf K}(x, y)=\frac{3\theta}{F}(x, y)+\kappa$, then $\widetilde{F}$ is of weakly isotropic flag curvature with $\widetilde{\bf K}(x, u)=\frac{3\tilde{\theta}}{\widetilde{F}}(x, u)+\tilde{\kappa}$. In this case,  $\theta =0$  and $\tilde{\theta}= - \rho _{x^m}u^{m}$, $\tilde{\kappa}= \kappa - \rho ^{2}+ 2 \rho _{x^m}(W^{m}+ V^{m})$. In particular, when $n \geq 3$, $\tilde{\theta}=0$ and $\tilde{\kappa}= \kappa - \rho ^{2}$ is a constant.
\item[{\rm (3)}] If $F$ is a weak Einstein metric with Ricci curvature
\be
Ric(x, y)=(n-1)\left\{\frac{3\theta}{F}+\kappa\right\} F^{2}(x, y),\label{WEKropina}
\ee
then $\widetilde{F}$ is a weak Einstein metric with Ricci curvature
\be
\widetilde{Ric}(x, u)=(n-1)\left\{\frac{3\widetilde{\theta}}{\widetilde{F}(x, u)}+\widetilde{\kappa}\right\}\widetilde{F}^{2}(x, u). \label{WERan}
\ee
Here,
\beqn
&& \tilde{\theta}= - \rho_{x^m}u^{m}, \\
&& \tilde{\kappa}= \kappa + \frac{3}{2}\theta _{m}W^{m}-\rho^{2}+2\rho_{x^m} (W^{m}+ V^{m})
\eeqn
and  $\kappa = \mu - \frac{3}{2}\theta _{i}W^{i}$ is a nonnegative scalar function on $M$ . Besides, here and after,
\[
A_{x}:=\left\{ y\in T_{x}M \ | \ h(y, W_{x})>0 \right\}=\left\{ y\in T_{x}M \ | \ \beta = b_{i}(x) y^{i} >0 \right\}.
\]
\een
\end{thm}

\vskip 2mm

Similarly, when $F(x, -V_{x})=1$, we have the following theorem.

\begin{thm}\label{Kropina's curvature} \ Let $F=\frac{\alpha^2}{\beta}$ be a Kropina metric on an $n (\geq 2)$ dimensional manifold $M$ expressed by (\ref{eqa3}) with navigation data $(h, W)$. Assume that $V$ is a Killing vector field on $(M, F)$ with $-V_{x}\in A_{x}$ and $F(x, -V_{x})=1$. Let $\widetilde{F}=\widetilde{F}(x, y)$ be the Kropina metric defined from $(F, V)$ by (\ref{KropinaNaviga}).  Then we have the following.
\ben
\item[{\rm (1)}] If $F$ is of isotropic $S$-curvature with ${\bf S}(x, y)=(n+1)cF(x, y)$, then $\widetilde{F}$ is of isotropic $S$-curvature with $\widetilde{\bf S}(x, u)=(n+1)\tilde{c}\widetilde {F}(x, u)$. Actually, in this case, $c= \tilde{c}=0$. Here and below, $u:=y+F(x, y)V=y+\widetilde{F}(x, u)V$.
\item[{\rm (2)}] If $F$ is of weakly isotropic flag curvature with ${\bf K}(x, y)=\frac{3\theta}{F}(x, y)+ \kappa$, then $\widetilde{F}$ is of weakly isotropic flag curvature with $\widetilde{\bf K}(x, u)=\frac{3\tilde{\theta}}{\widetilde{F}}(x, u)+\tilde{\kappa}$. In this case, $\theta = \tilde{\theta}=0$ and $\tilde{\kappa} = \kappa \geq 0$.
\item[{\rm (3)}]  If $F$ is a weak Einstein metric with Ricci curvature
\be
      Ric(x, y)=(n-1)\left\{\frac{3\theta}{F}+\kappa \right\}F^{2}(x, y),\label{WEKropina2}
\ee
then $\tilde{F}$ is also a weak Einstein metric with Ricci curvature
\[
\widetilde{Ric}(x, u)=(n-1)\left\{\frac{3\tilde{\theta}}{\widetilde{F}}(x, u) +\tilde{\kappa}\right\}\widetilde{F}^{2}(x, u). \label{NWEKropina}
\]
\een
where
\[
 \tilde{\kappa}(x)= \kappa(x)+\frac{3}{2}\left\{(\theta _{m}-\tilde{\theta}_{m})W^{m}- \tilde{\theta}_{m}V^{m}\right\}.
\]
\end{thm}

\vskip 2mm

The paper is organized as follows. In Section \ref{sec2}, we give some definitions and curvature properties of Randers metrics which are necessary for the present paper. Some useful results on Kropina metrics are given in Section \ref{sec3}. In particular, we characterize weak Einstein Kropina metrics via navigation data in this section. Furthermore, we prove that the solution of Zermelo navigation problem on any Kropina manifold $(M, F)$ with a vector field $V$ satisfying $F(x, -V_{x})\leq 1$ must be either a Randers metric or a Kropina metric in Section \ref{sec4}. Finally, we give the proofs of the main theorems in this paper in Section \ref{sec5}.

\section{Preliminaries}\label{sec2}

Let $F$ be a Finsler metric on an $n$-dimensional smooth manifold $M$ and $(x, y)=(x^{i}, y^{i})$ the local coordinates on the tangent bundle $TM$. Let
\[
G^{i}=\frac{1}{4}g^{il}\left\{[F^{2}]_{x^{m}y^{l}}y^{m}-[F^{2}]_{x^{l}}\right\}
\]
be the spray coefficients of $F$. For any $x\in M$ and $y\in T_{x}M\setminus \{0\}$, the Riemann curvature ${\bf R}= R^{i}_{\ k}\frac{\partial}{\partial x^{i}}\otimes dx^{k}$ is defined by
\[
R^{i}_{\ k}=2\frac{\partial G^{i}}{\partial x^{k}}-y^{m}\frac{\partial^{2}G^{i}}{\partial x^{m}\partial y^{k}}+2G^{m}\frac{\partial^{2}G^{i}}{\partial y^{m}\partial y^{k}}-\frac{\partial G^{i}}{\partial y^{m}}\frac{\partial G^{m}}{\partial y^{k}}.
\]
For Finsler manifold $(M, F)$, the flag curvature ${\bf K}={\bf K} (P, y)$  is a function of ``flag" $P \subset T_{x}M $ and ``flagpole" $y\in T_{x}M $ at $x \in M$ with $y\in P $. The flag curvature in Finsler geometry is a natural extension of sectional curvature in Riemannian geometry and is the most important Riemannian geometric quantity in Finsler geometry. A Finsler metric $F$ is said to be of scalar flag curvature if
${\bf K}={\bf K}(x,y)$ is independent of the flag $P$. In particalar, we say that Finsler metric $F$ is of weakly isotropic flag curvature if
\be
{\bf K}=\frac{3\theta}{F}+\sigma,
\ee
where $\sigma=\sigma (x)$ is a scalar function and $\theta = \theta _{i}y^{i}$ is a 1-form on $M$.  It is known that $F$ is of scalar flag curvature if and only if, in a standard local coordinate system,
\[
R^{i}_{\ k}= {\bf K} (x, y)\{F^{2}\delta^{i}_{k}-FF_{y^{k}}y^{i}\}.
\]
The Ricci curvature of $F$ is defined by $Ric =R^{m}_{\ m}$. $F$ is called a weak Einstein metric if
\be
Ric = (n-1)\left\{\frac{3\theta}{F} + \kappa \right\}F^{2},
\ee
where $\theta$ is a 1-form and $\kappa =\kappa (x)$ is a scalar function on $M$. When $\theta =0$, $F$ is called an Einstein metric.

$S$-curvature is a very  important non-Riemannian quantity in Finsler geometry, which was introduced by Z. Shen when he extended Bishop-Gromov volume comparison theorem to Finsler geometry (see \cite{ZS1}). For a volume form $dV=\sigma _{F}(x)dx^{1} \cdots dx^{n}$ on $M$, S-curvature is given by
\[
{\bf S}(x, y)=\frac{\partial G^{m}}{\partial y^{m}}-y^{m}\frac{\partial \ ln\sigma_{F}}{\partial x^{m}}.
\]
The $S$-curvature measures the average rate of change of tangent space $(T_{x}M, F_{x})$ in direction $y\in T_{x}M$. If ${\bf S} =(n+1)\{cF+\eta\}$, where $c=c(x)$ is a scalar function  and $\eta$ is a closed 1-form on $M$,   $F$ is said to be of almost isotropic $S$-curvature.  If $\eta =0$, that is, ${\bf S}=(n+1)cF$, we say that $F$ has isotropic $S$-curvature (see \cite{ChernShen}).

For Randers metrics, we have the following important proposition which is necessary for the proofs of our main theorems.

\begin{prop}{\rm (\cite{ChernShen})}\label{Chern-shen-Prop}\
Let $F=\alpha+\beta$ be a Randers metric on a manifold $M$ of dimension $n$ expressed by (\ref{RhW}) with navigation data $(h, W)$. Then $F$ has isotropic $S$-curvature, ${\bf S}=(n+1)c(x)F$, if and only if $W$ satisfies
\beq
W_{i|j}+W_{j|i} = - 4c(x)h_{ij}
\eeq
where $c(x)$ is a scalar function on $M$ and $``|"$ denotes the covariant derivative with respect to the Levi-Civita connection of  $h$.
\end{prop}

It is a long standing  important problem in Finsler geometry to classify Finsler metrics of scalar flag curvature. Bao-Robles-Shen have classified completely strongly convex Randers metrics of constant flag curvature(\cite{BRS}). Further, the first author and X. Mo, Z. Shen proved that, for a Finsler metric $F$ of scalar flag curvature, if it is of isotropic S-curvature, then $F$ must be of weakly isotropic flag curvature (\cite{CSM}). In general, the converse of this conclusion does not necessarily hold except for Randers metrics and Kropina metrics (\cite{XZ}\cite{XQ2}). Base on the result, the first author and Z. Shen have classified completely the Randers metrics of weakly isotropic flag curvature by navigation data.

\begin{prop}{\rm (\cite{XZ})}\label{flag curvature}
Let $F=\alpha+\beta$ be a Randers metric on a manifold $M$ of dimension $n$ expressed by (\ref{RhW}) with navigation data $(h, W)$. Then $F$ is of weakly isotropic flag curvature
\be
{\bf K} =\frac{3c_{x^m}y^{m}}{F}+\sigma \label{weakflag}
\ee
if and only if $F$ has isotropic $S$-curvature  ${\bf S}=(n+1)cF$ and $h$ is of isotropic sectional curvature ${\bf K}_{h}=\mu (x)$,
where $c(x)$ and $\sigma (x)$ are scalar functions on $M$ and $\mu = \sigma(x) +c^{2}+2c_{x^m}W^{m}$.
\end{prop}

Besides, the first author and Z. Shen have also characterized the weak Einstein Randers metrics.

\begin{prop}{\rm (\cite{XZ})}\label{Ricci curvature}
Let $F=\alpha+\beta$ be a Randers metric on a manifold $M$ of  dimension $n$  expressed by (\ref{RhW}) with navigation data $(h, W)$. Assume that $F$ is of isotropic $S$-curvature, ${\bf S} =(n+1)cF$. Then $F$ is weak Einstein metric with
\be
Ric=(n-1)\left\{\frac{3c_{x^m}y^{m}}{F}+\sigma \right\}F^{2},
\ee
if and only if $h$ is an Einstein metric with
\be
{}^{h}Ric =(n-1)\mu h^{2},
\ee
where $c= c(x)$ and $\sigma =\sigma (x)$ are scalar functions on manifold $M$ and $\mu = \sigma(x) +c^{2} +2c_{x^m}W^{m}$, ${}^{h}Ric$ is Ricci curvature of Riemannian metric $h$.
\end{prop}

\section{Some useful results on Kropina metrics}\label{sec3}

In this section, we will introduce some useful results on Kropina metrics which are important for the discussions below.

\begin{lem}{\rm (\cite{XQ2})}\label{S-curvature}
For a Kropina metric $F=\frac{\alpha ^{2}}{\beta}$ on an $n$-dimensional manifold $M$, the following are equivalent.
\ben
\item[{\rm (a)}] $F$ has an isotropic $S$-curvature, that is, ${\bf S} =(n+1)cF$;
\item[{\rm (b)}] $r_{00}=\sigma \alpha^{2}$;
\item[{\rm (c)}] ${\bf S} =0$;
\item[{\rm (d)}] $\beta$ is a conformal form with respect to $\alpha$,
\een
where $c=c(x)$ and $\sigma = \sigma (x)$ are functions on $M$ and $r_{00}:=r_{ij}y^{i}y^{j}$, $r_{ij}:=\frac{1}{2}(b_{i;j}+b_{j;i})$ and $`` ; "$ denotes the covariant derivative with respect to the Levi-Civita connection of $\alpha$.
\end{lem}

Lemma \ref{S-curvature} tells us an important fact that a Kropina metric $F$ is of isotropic S-curvature if and only if its S-curvature vanishes, ${\bf S}=0$.  On the other hand,  X. Zhang and Y. Shen got the following result (see Lemma 4.1, \cite{ZY}).
\begin{lem}{\rm (\cite{ZY})}\label{RSW} \ For a Kropina metric $F=\frac{\alpha ^2}{\beta}$ with navigation data $(h, W)$,
$r_{00}= \sigma(x)\alpha^{2}$ is equivalent to the following
\be
W_{i|j}+W_{j|i}=0. \label{Kr00}
\ee
Here, $``|"$ denotes the covariant derivative with respect to the Levi-Civita connection of  $h$.
\end{lem}

Consequently, by Lemma \ref{S-curvature} and Lemma \ref{RSW},  we get the following proposition.

\begin{prop}\label{corA}\
Let $F=\frac{\alpha^2}{\beta}$ be a Kropina metric on manifold $M$ expressed by (\ref{eqa3}) via navigation data $(h, W)$. Then $F$ is of isotropic $S$-curvature if and only if
\be
W_{i|j}+W_{j|i}=0, \label{c0}
\ee
that is, $W$ is a Killing vector field with respect to $h$. In this case, ${\bf S} =0$.
\end{prop}

In 2013, R. Yoshikawa and  S. V. Sabau classified the Kropina metrics of constant flag curvature by navigation data $(h, W))$ and shown that, up to local isometry, there are only two model spaces of them: the Euclidean space and the odd-dimensional spheres (\cite{RS}). More generally, Q. Xia characterized the Kropina metrics of weakly isotropic flag curvature and proved the following result.

\begin{lem}{\rm (\cite{XQ2})}\label{KFC}
Let $F=\frac{\alpha^{2}}{\beta}$ be a Kropina metric on an $n(\geq2)$-dimensional manifold $M$ with the navigation data $(h, W)$. Then $F$ is of weakly isotropic flag curvature ${\bf K} =\frac{3\theta}{F}+\kappa$ if and only if the sectional curvature of $h$ is a nonnegative scalar function $\mu(x)$ and $W$ is a Killing vector field with respect to $h$. In this case, ${\bf K} = \kappa = \mu \geq 0$ and $\theta=0$.
\end{lem}

In \cite{ZY}, X. Zhang and Y. Shen studied Enistein Kropina metrics and proved that a Kropina metric $F$ on an $n(\geq2)$-dimensional manifold $M$ with the navigation data $(h, W)$  is an Einstein metric with $Ric = (n-1)\sigma F^{2}$ if and only if $h$ is an Einstein metric with
${}^{h}Ric = (n-1)\mu h^{2}$ and $W$ is a Killing vector field with respect to $h$. In this case, $\sigma = \mu \geq 0$. Moreover, $\sigma$ is a constant for $n \geq 3$. Further, it is not difficult to prove the following lemma which characterizes weak Einstein Kropina metrics and was first mentioned in  Remark 5.2 of \cite{XQ2}.

\begin{lem} \label{Ricci}
Let $F=\frac{\alpha^{2}}{\beta}$ be a Kropina metric on an $n(\geq2)$-dimensional manifold $M$ and $b:=\|\beta\|_{\alpha}$ be a constant. Then $F$ is a weak Einstein metric with
\[
Ric=(n-1)\left\{\frac{3\theta}{F}+\kappa\right\}F^{2}
\]
if and only if $\alpha$ is an Einstein metric with Ricci scalar $\mu =\mu(x)$ and $\beta$ is a Killing 1-form. In this case, $\kappa=\frac{1}{4}(\mu b^{2}-3\theta_{i}b^{i})\geq 0$.
\end{lem}

Based on Lemma \ref{Ricci}, we can characterize weak Einstein Kropina metrics via navigation data $(h, W)$ from (\ref{eqa3}).

\begin{thm}\label{Einstein's prop}
Let $F=\frac{\alpha^2}{\beta}$ be a Kropina metric on an $n(\geq2)$-dimensional manifold $M$ expressed by (\ref{eqa3}) with navigation data ($h$, $W$). Then $F$ is a weak Einstein metric with
\be
Ric=(n-1)\left\{\frac{3\theta}{F}+\kappa\right\}F^{2} \label{weakEK}
\ee
if and only if $h$ is an Einstein metric with ${}^{h}Ric=(n-1)\mu (x) h^{2}$ and $W$ is a Killing vector field with respect to $h$, where $\theta = \theta _{i}y^{i}$ is a 1-form and $\kappa$ is a nonnegative scalar function on $M$. In this case, $\kappa=\mu-\frac{3}{2}\theta_{i}W^{i}\geq 0$.
\end{thm}
\noindent {\bf Proof.} \ From $F=\frac{h^{2}}{2W_{0}}$, let $\bar{F}:=\frac{h^2}{W_{0}}$. Then $F=\frac{1}{2}\bar{F}$ and $Ric= \overline{Ric}$, where $\overline{Ric}$ denotes the Ricci curvature of $\bar{F}$. In this case, $F$ is a weak Einstein metric with (\ref{weakEK}) if and only if $\bar{F}$ is a weak Einstein metric with
\be
\overline{Ric}=(n-1)\left\{\frac{3\bar{\theta}}{\bar{F}}+\bar{\kappa}\right\}\bar{F}^{2}, \label{weakEKbar}
\ee
where $\bar{\theta}=\frac{1}{2}\theta$ and $\bar{\kappa}=\frac{1}{4}\kappa$.

On the other hand, noticing that $\|W_{0}\|_{h}=1$, by Lemma \ref{Ricci}, $\bar{F}$ is a weak Einstein metric with (\ref{weakEKbar}) if and only if $h$ is an Einstein metric with ${}^{h}Ric=(n-1)\mu (x) h^{2}$ and $W_{0}$ is a Killing 1-form, and $\bar{\kappa}=\frac{1}{4}\{\mu (x)\|W_{0}\|^{2}_{h}-3 \bar{\theta}_{i}W^{i}\}=\frac{1}{4}\{\mu (x) -\frac{3}{2}\theta _{i}W^{i}\}\geq 0$. Then we know that $W$ is a Killing vector field and $\kappa=\mu-\frac{3}{2}\theta_{i}W^{i}\geq 0$. This completes the proof of the theorem. \qed

\vskip 2mm

By Proposition \ref{corA}, Lemma \ref{KFC} and Theorem \ref{Einstein's prop}, it is easy to get the following proposition.

\begin{prop}\label{wifEinS}
Every  Kropina metric of weakly isotropic flag curvature or weak Einstein Kropina metric on an $n (\geq2)$-dimensional manifold has vanishing S-curvature, ${\bf S}=0.$
\end{prop}

A vector field $V$ on a Finsler manifold $(M,F)$ is called a conformal vector field
with  conformal factor $\rho = \rho (x)$ if the one-parameter transformation group $\{\varphi_{t}\}$ generated by $V$ is a conformal transformation group, that is,
\be
F\left(\varphi_{t}(x), (\varphi_{t})_{*}(y)\right)=e^{2 \sigma_{t}(x)} F(x, y), \ \ \forall x \in M, \ y \in T_{x} M, \label{confvector}
\ee
where $\sigma_{t}(x) :=\int_{0}^{t} \rho\left(\varphi_{s}(x)\right) ds$. In this case, it is easy to see that $\rho (x)=\frac{d \sigma_{\mathrm{t}}(x)}{d t} |_{t=0}$ and $\sigma_{0}(x)=0$. In particular, $V$ is called a  homothetic vector field on $M$ if $\rho$ is constant. $V$ is called a  Killing vector field if $\rho=0$.

\begin{lem}{\rm (\cite{CLY})}\label{lem2.1}
Let $F=\frac{\alpha^2}{\beta}$ be a Kropina metric on manifold $M$ given by (\ref{eqa3}) with navigation data $(h, W)$. Then a vector field $V$ on $(M,F)$ is a conformal vector field with conformal factor $\rho=\rho(x)$ if and only if $V$ satisfies the following conditions:
\beq
&&V_{i|j}+V_{j|i}=4\rho h_{ij},\label{c1}\\
&&V^i W_{j|i}+W^i V_{i|j} =2\rho W_{j}, \label{c2}
\eeq
where we use $h_{ij}$ to raise and lower the indices of $V$ and $W$ and $`` \ | \ "$ denotes the covariant derivative with respect to the Levi-Civita connection of  $h$.
\end{lem}

From Lemma \ref{lem2.1}, if $V$ is a conformal vector field with conformal factor $\rho = \rho(x)$ on $(M,F=\frac{\alpha ^2}{\beta} )$ , then $V$ is a conformal vector field with conformal factor $\rho = \rho (x)$ on $(M, h)$.

\section{An important fundamental theorem}\label{sec4}

In this section, we will give an important theorem, which will provide an important basis for Theorem \ref{Kropina navigation} and Theorem \ref{Kropina's curvature}

By the proof of Lemma 4.1 in \cite{SX}, we can actually get the following useful lemma. 

\begin{lem} \label{lemma3.1} \
Let $ \Phi (x, y)$ be a Finsler metric on a manifold $M$ and $W$ be a vector field on $M$ with $\Phi (x, -W_{x}) \leq 1$. Suppose $F(x, y)$ is a Finsler metric defined from navigation data $(\Phi, W)$ by (\ref{navigation}) and $V$ is a vector field on $M$ with $F(x, -V_{x}) \leq 1$. Then the Finsler metric $\widetilde{F}(x, y)$ defined from navigation data $(F, V)$ by
\be
F\left(x, \frac{y}{\widetilde{F}(x,y)}- V_{x}\right)=1 \label{navigaFV}
\ee
satisfies the following identity,
\be
\widetilde{F}(x, u)= \Phi (x, u-\widetilde{F}(x, u)(W+V)), \label{navigawideFVW}
\ee
where $u=y+ F(x, y)V$. That is, $\widetilde{F}$ is just the solution of navigation problem with navigation data $(\Phi , W+V)$. In this case, $\widetilde{F}(x, u)= F(x,y)$.
\end{lem}

According to Lemma \ref{lemma3.1}, we can determine the solution of Zermelo navigation problem on a Kropina manifold $(M, F)$ with navigation data $(h, W)$. In this case, $\|W\|_{h}=1$. Hence, at each point $x\in M$, the indicatrix $S_{F}=S_{h}+\{ W\}$ passes thorough the origin $x$ of $T_{x}M$.

\begin{thm}\label{lem3}\
Let $F=\frac{\alpha^2}{\beta}$ be a Kropina metric on manifold $M$ given by (\ref{eqa3}) with navigation data ($h$, $W$). Suppose $V$ is a vector field on $M$ with  $- V_{x}\in A_{x}$ and $F(x, -V_{x})\leq 1$. Let $\widetilde{F}$ be the solution of (\ref{navigaFV}) with navigation data $(F, V)$. Then we have the following.
\ben
\item[{\rm (1)}] If $F(x, -V_{x})<1$, then $\widetilde{F}$ is a Randers metric.
\item[{\rm (2)}] If $F(x, -V_{x})=1$, then $\widetilde{F}$ is a  Kropina metric.
\een
\end{thm}
\noindent {\bf Proof.} \ By the assumption, we have
\[
F=\frac{\alpha^2}{\beta} = \frac{h^{2}}{2W_{0}},
\]
where $h^{2} =\frac{4}{b^2}\alpha ^2$, $W_{0}=\frac{2}{b^{2}}\beta$, $b:=\|\beta\|_{\alpha}$ and $\|W\|_{h}=1$. In this case, the assumption that $- V_{x}\in A_{x}$ implies the following
\be
F(x, -V_{x})=\frac{\|V_{x}\|_{h}^{2}}{-2W_{i}V^{i}}>0. \label{lengthmV}
\ee
Furthermore, by Lemma \ref{lemma3.1} and the assumption, $\widetilde{F}$ is just the solution of the navigation problem with navigation data $(h, W+V)$.

{\it Case I}.  $F(x, -V)<1$. In this case, by (\ref{lengthmV}),
\[
\|V_{x}\|_{h}^{2}+2W_{i}V^{i}<0.
\]
 By $\|W\|_{h}=1$, we get
\[
\|V_{x}\|_{h}^{2}+2W_{i}V^{i}+\|W\|_{h}^{2}<1,
\]
which implies that $\|V+W\|_{h}<1$. Hence,  $\widetilde{F}$ is a  Randers metric with navigation data $(h, V + W)$.

{\it Case II}.  $F(x, -V)= 1$. In this case, by (\ref{lengthmV}),
\[
\|V_{x}\|_{h}^{2}+2W_{i}V^{i}=0.
\]
By $\|W\|_{h}=1$, we get
\[
\|V_{x}\|_{h}^{2}+2W_{i}V^{i}+\|W\|_{h}^{2} =1,
\]
which implies that $\|V+W\|_{h} = 1$. Hence,  $\widetilde{F}$ is a new Kropina metric with navigation data $(h, V + W)$.
\qed

\vskip 2mm

The following corollary is obvious.
\begin{cor}
Let $F=\frac{\alpha^2}{\beta}$ be a Kropina metric on manifold $M$ given by (\ref{eqa3}) with navigation data ($h$, $W$). Suppose $V$ is a vector field on $M$ with  $- V_{x}\in A_{x}$ and $F(x, -V_{x})\leq 1$. Let $\widetilde{F}$ be the solution of (\ref{navigaFV}) with navigation data $(F, V)$. Then we have the following.
\ben
\item[{\rm (1)}] If $F(x, -V_{x})<1$, then $\widetilde{F}$ is a Randers metric given by
\be
\widetilde{F}=\frac{\sqrt{\tilde{\lambda} {h}^{2}+\widetilde{W}_{0}^{2}}}{\tilde{\lambda}}+\frac{\widetilde{W_{0}}}{\tilde{\lambda}}, \label{d2}
\ee
where $\widetilde{W_{0}}:=V_{0}+W_{0}$ and $\tilde{\lambda}:=-(\|V_{x}\|_{h}^{2}+2 W_{i}V^{i})>0$.

\item[{\rm (2)}] If $F(x, -V_{x})=1$, then $\widetilde{F}$ is a new Kropina metric given by
\be
\widetilde{F}=\frac{h^{2}}{2\widetilde{W}_{0}}, \label{d23}
\ee
where $\widetilde{W_{0}}:=V_{0}+W_{0}$. In this case, $\|V+W\|_{h}=1$.
\een
\end{cor}

\begin{rem} In Theorem \ref{lem3}, the condition that $- V_{x}\in A_{x}$ means that, at each point $x \in M$,  the parallel shift $S_{F}+\{V\}$ of the indicatrix of $F$ contains the origin $x$ of $T_{x}M$ when $F(x, - V_{x})<1$, while  $S_{F}+\{V\}$ passes through the origin $x$ of $T_{x}M$ when $F(x, -V_{x})=1$. In fact, when $F(x, -V_{x})=1$, the endpoint of $-V_{x}$ on $S_{F}$ will coincide with the origin $x$ after parallel shift along $V_{x}$.
\end{rem}

\section{Proofs of Theorem \ref{Kropina navigation} and Theorem \ref{Kropina's curvature}}\label{sec5}

In this section, we will give the proofs of Theorem \ref{Kropina navigation} and Theorem \ref{Kropina's curvature}.

\vskip 2mm

\noindent
{\it Proof of Theorem \ref{Kropina navigation}.} \ Firstly, we should remember that $\widetilde{F}$ is just the solution of the navigation problem with navigation data $(h, W+V)$ satisfying $\|W+V \|_{h}<1$ by Theorem \ref{lem3}.

By the assumption, $V$ is a conformal vector field of $F$ with conformal factor $\rho(x)$. Then, by Lemma \ref{lem2.1}, we have
\be
V_{i|j}+V_{j|i}=4\rho h_{ij}, \label{e1}
\ee
which implies that $V$ must be a conformal vector field of $h$ with conformal factor $\rho(x)$.

\vskip 1mm

(1)\ Since $F$ is of isotropic $S$-curvature,  ${\bf S} =(n+1)c(x)F$, by Proposition \ref{corA},
\be
W_{i|j}+W_{j|i}=0, \label{thm1.1S}
\ee
that is, $W$ is a Killing vector field with respect to $h$. In this case, ${\bf S} =0$, that is, $c=0$.

Adding (\ref{e1}) and (\ref{thm1.1S}) yields
\be
(V_{i}+W_{i})_{|j}+(V_{j}+W_{j})_{|i}=4\rho h_{ij}. \label{randersS}
\ee
By Proposition \ref{Chern-shen-Prop} and (\ref{randersS}), we know that $\widetilde{F}$ has isotropic S-curvature with $\widetilde{\bf S}=(n+1)(- \rho)\widetilde{F}$.

\vskip 1mm

 (2)\ Since $F$ is of weakly isotropic flag curvature,  ${\bf K}(x, y)=\frac{3\theta}{F}+ \kappa$, by Lemma \ref{KFC}, the sectional curvature of $h$ is a nonnegative scalar function, ${\bf K}_{h}=\mu(x)\geq 0$ and $W$ is a Killing vector field on $M$. In this case, $\theta =0$ and  ${\bf K}(x, y)= \kappa = \mu $.

From Proposition \ref{wifEinS}, $F$ has vanishing S-curvature. Then, by (1), we know that  $\widetilde{F}$ has isotropic S-curvature with $\widetilde{\bf S}=(n+1)(- \rho)\widetilde{F}$. Thus, by Proposition \ref{flag curvature}, we can conclude that $\widetilde{F}$ is of weakly isotropic flag curvature with
 \[
 \widetilde{\bf K}= \frac{3\tilde{\theta}}{\widetilde{F}}(x, y)+\tilde{\kappa},
 \]
where $\tilde{\theta}=- \rho _{x^m}y^{m}$ and $\tilde{\kappa}=\kappa - \rho ^{2}+2 \rho _{x^{m}}(W^{m}+ V^{m})$. In particular, when $n \geq 3$, by Theorem 1.3 in \cite{CLY}, $V$ is a homothetic vector field on $M$, which means that $\rho$ is a constant. Hence, in this case, $\tilde{\theta}=0$ and $\tilde{\kappa} = \kappa - \rho ^{2}$ is a constant.

\vskip 1mm

 (3)\ Because $F$ is a weak Einstein metric with (\ref{WEKropina}), from Theorem \ref{Einstein's prop}, $h$ is an Einstein metric  with Ricci curvature ${}^{h}Ric=(n-1)\mu (x) h^{2}$ and $W$ is a Killing vector field. In this case, $\kappa=\mu-\frac{3}{2}\theta_{i}W^{i}\geq 0$ and $F$ has isotropic $S$-curvature from Proposition \ref{corA}.

 By (1) again, we konw that $\widetilde{F}$ is also of isotropic $S$-curvature with $\widetilde{\bf S}=(n+1)(- \rho)\widetilde{F}$. Then, from Proposition \ref{Ricci curvature}, we can assert that $\widetilde{F}$ is a weak Einstein metric with
\[
\widetilde{Ric}=(n-1)\left\{\frac{3\widetilde{\theta}}{\widetilde{F}}+\widetilde{\kappa}\right\}\widetilde{F}^{2},
\]
where $\tilde{\theta}= - \rho_{x^m}y^{m}$, $\widetilde{\kappa}:=\mu- \rho^{2}+ 2 \rho_{x^m}(W^{m}+ V^{m})= \kappa+\frac{3}{2}\theta_{m}W^{m}-\rho^{2}+ 2 \rho_{x^m}(W^{m}+V^{m})$.

In sum, considering the corresponding relationship between the points $(x, y)$ and $(x, u)$ established in Lemma \ref{lemma3.1}, we obtain Theorem \ref{Kropina navigation}.  \qed

\vskip 2mm

The following is the proof of Theorem \ref{Kropina's curvature}.
\vskip 2mm
\noindent {\it Proof of Theorem \ref{Kropina's curvature}.}  Similar to the proof of Theorem \ref{Kropina navigation},  $\widetilde{F}$ is just the solution of the navigation problem with navigation data $(h, W+V)$ satisfying $\|W+V \|_{h} = 1$ by Theorem \ref{lem3}.

By the assumption, $V$ is a Killing vector field of $F$. Then, by Lemma \ref{lem2.1}, $V$ must be a Killing vector field of $h$, that is,
\be
V_{i|j}+V_{j|i}=0. \label{f1}
\ee

\vskip 1mm

(1)\ If $F$ is of isotropic $S$-curvature with ${\bf S}(x, y)=(n+1)cF(x, y)$, by Proposition \ref{corA}, $W$ is a Killing vector field with respect to $h$,
\be
W_{i|j}+W_{j|i}=0. \label{KropinaS-cur}
\ee
In this case, ${\bf S}=0$,  that is, $c=0$.

Adding (\ref{f1}) and (\ref{KropinaS-cur}) yields
\[
(V_{i}+W_{i})_{|j}+(V_{j}+W_{j})_{|i}=0,
\]
 which means that $V+W$ is still a Killing vector field of $h$. By Proposition \ref{corA} again, we know that $\widetilde{F}$ has vanishing S-curvature, ${\widetilde{\bf S}}=0$.

\vskip 1mm

 (2)\ If $F$ is of weakly isotropic flag curvature,  ${\bf K}(x, y)=\frac{3\theta}{F}+ \kappa$, by  Lemma \ref{KFC}, the sectional curvature of $h$ is a nonnegative scalar function ${\bf K}_{h}= \mu(x)$ and $W$ is a Killing vector field with respect to $h$. In this case, $\theta =0$ and  ${\bf K}(x, y)= \kappa = \mu \geq 0$.

From (1), $W+V$ is still  a Killing vector field of $h$. Then, by  Lemma \ref{KFC} again, we conclude that $\widetilde{F}$ is of weakly isotropic flag curvature
$\widetilde{\bf K}(x, y)=\frac{3\tilde{\theta}}{\widetilde{F}(x, y)}+\tilde{\kappa}$ and $\tilde{\theta}=0$,   $\widetilde{\bf K}(x, y)=\tilde{\kappa}= \kappa =\mu(x) \geq 0$.

\vskip 1mm

(3)\ Because $F$ is a weak Einstein metric with (\ref{WEKropina2}), $h$ is an Einstein metric with   ${}^{h}Ric=(n-1)\mu (x) h^{2}$  and $W$ is a Killing vector field from Theorem \ref{Einstein's prop}.  In this case, $\kappa=\mu-\frac{3}{2}\theta_{i}W^{i}\geq 0$ and $F$ has isotropic $S$-curvature from Proposition \ref{corA}.

Further, it is  easy to see that $W+V$ is still a Killing vector field of $h$. Hence, by Theorem \ref{Einstein's prop} again, $\widetilde{F}$ is a weak Einstein metric with
\[
\widetilde{Ric}(x, y)=(n-1)\left\{\frac{3\widetilde{\theta}}{\widetilde{F}(x, y)}+\widetilde{\kappa}\right\}\widetilde{F}^{2}(x, y).
\]
Here $\widetilde{\theta}$ is a 1-form and $\widetilde{\kappa}$ is a nonnegative scalar function on $M$. In this case, $\tilde{\kappa}$:=$\mu-\frac{3}{2}\tilde{\theta}_{m}(W^{m}+V^{m})=\kappa+\frac{3}{2}\left\{(\theta_{m}-\widetilde{\theta}_{m})W^{m}-\tilde{\theta}_{m}V^{m})\right\}$.

Also, because of the corresponding relationship between the points $(x, y)$ and $(x, u)$ established in Lemma \ref{lemma3.1}, we obtain Theorem \ref{Kropina's curvature}. \qed

\vskip 8mm

\noindent
Xinyue Cheng \\
School of Mathematical Sciences \\
Chongqing Normal University \\
Chongqing  401331,  P. R. of China  \\
E-mail: chengxy@cqnu.edu.cn

\vskip 4mm

\noindent
Qiuhong Qu \\
School of Mathematical Sciences \\
Chongqing Normal University \\
Chongqing  401331,  P. R. of China  \\
E-mail: 465206186@qq.com

\vskip 4mm

\noindent
Suiyun Xu \\
School of Sciences\\
Chongqing University of Technology\\
Chongqing, 400054, P. R. China\\
E-mail: 764536895@qq.com

\end{document}